

\input amstex

\documentstyle{amsppt}

\loadbold

\magnification=\magstep1

\pageheight{9.0truein}
\pagewidth{6.5truein}


\def\k{\overline{k}}
\def\Z{\Bbb{Z}}
\def\Q{\Bbb{Q}}

\def\nxn{n{\times}n}
\def\ker{\operatorname{ker}}
\def\M{\Bbb{M}}
\def\Rel{\operatorname{Rel}}
\def\sgn{\operatorname{sgn}}
\def\trace{\operatorname{trace}}
\def\x{\bold{x}}
\def\y{\bold{y}}
\def\z{\bold{z}}

\def\Art{1}
\def\BouBou{2}
\def\CoxLitOsh{3}
\def\Forone{4}
\def\Fortwo{5}
\def\Fri{6}
\def\Lab{7}
\def\Let{8}
\def\Lin{9}
\def\McCRob{10}
\def\Pap{11}
\def\PleSou{12}
\def\Pro{13}
\def\Raz{14}
\def\Row{15}

\topmatter

\title Constructing Irreducible Representations of Finitely Presented
Algebras \endtitle

\rightheadtext{Constructing Irreducible Representations}

\author Edward S. Letzter \endauthor

\abstract We describe an algorithmic test, using the ``standard
polynomial identity'' (and elementary computational commutative
algebra), for determining whether or not a finitely presented
associative algebra has an irreducible $n$-dimensional
representation. When $n$-dimensional irreducible representations do
exist, our proposed procedure can (in principle) produce explicit
constructions. \endabstract

\address Department of Mathematics, Temple University, Philadelphia,
PA 19122 \endaddress

\email letzter\@math.temple.edu\endemail

\thanks The author's research was supported in part by grants from the
National Science Foundation. \endthanks

\endtopmatter

\document

\head 1. Introduction \endhead

Our aim in this note is to suggest a general algorithmic approach to the
finite dimensional irreducible representations of finitely presented algebras,
combining well-known methods from both noncommutative ring theory and
computational commutative algebra. (There have been numerous previous studies,
from an algorithmic perspective, on matrix representations of finitely
presented groups and algebras; see, e.g., \cite{\Lab, \Lin, \PleSou} for
analyses that -- as in our study below -- do not place additional technical
conditions on the groups or algebras involved.)

\subhead 1.1 \endsubhead To briefly describe the content of this note, assume
that $k$ is a computable field, and that $\k$ denotes the algebraic closure of
$k$. Suppose further that $n$ is a fixed positive integer, that $M_n(\k)$ is
the algebra of $\nxn$ matrices over $\k$, and that $R$ is a finitely presented
$k$-algebra. We will always use the expression {\sl $n$-dimensional
representation of $R$\/} to mean a $k$-algebra homomorphism $\rho\colon R
\rightarrow M_n(\k)$, and we will say that $\rho$ is {\sl irreducible\/} when
$M_n(\k)$ is $\k$-linearly spanned by $\rho(R)$ (cf\., e.g., \cite{\Art, \S
9}). Note that $\rho$ is irreducible if and only if $\rho\otimes 1 \colon R
\otimes_k \k \rightarrow M_n(\k)$ is surjective, if and only if $\rho \otimes
1$ is irreducible in the more common use of the term.

\subhead 1.2 \endsubhead Calculating over $k$, the procedure
described in this note always (in principle)

\itemitem{(a)} decides whether irreducible representations $R \rightarrow
M_n(\k)$ exist,

\itemitem{(b)} explicitly constructs an irreducible representation $R
\rightarrow M_n(\k)$ if at least one exists (assuming that $k[x]$ is
equipped with a factoring algorithm).  \par

\subhead 1.3 \endsubhead The finite dimensional irreducible
representations of finitely generated noncommutative algebras were
parametrized, up to equivalence (i.e., up to isomorphisms among the
corresponding modules), in famous work of Artin \cite{\Art}, Formanek
\cite{\Forone}, Procesi \cite{\Pro}, Razmyslov \cite{\Raz}, and others
(see, e.g., \cite{\Fortwo}, \cite{\McCRob}, or \cite{\Row}). The
algorithm we describe, however, does not distinguish among equivalence
classes of irreducible representations; it depends only on the
Amitsur-Levitzky Theorem (e.g., \cite{\McCRob, 13.3.3}) and recent
work of Pappacena \cite{\Pap}. In \cite{\Let} we present a procedure
that counts the number (possibly infinite) of equivalence classes of
irreducible representations, in characteristic zero.

\subhead 1.4 \endsubhead Examples are discussed in section 4. All of
the computational commutative algebra used in this note is elementary,
and the necessary background can be found in \cite{\CoxLitOsh}, for
example.

\subhead Acknowledgements \endsubhead I am grateful to Martin Lorenz,
Rekha Thomas, Zinovy Reichstein, and Bernd Sturmfels for reading
previous drafts of this note and for their very helpful comments. I am
also grateful to the referees, for suggesting several substantial
revisions that have improved both the content and exposition of this
note; see, in particular, (2.1), (2.2), and (2.6vii).

\head 2. Representations of Finitely Presented Algebras \endhead

While most of the material in this section is known, we provide a
complete treatment, for the reader's convenience.

\subhead 2.0 \endsubhead (i) Retaining the notation of (1.1), let $\M$
denote the affine space $(M_n(\k))^s$ of $s$-tuples of $\nxn$ matrices
over $\k$.

(ii) Let
$$B = k[x_{ij}(\ell) : 1 \leq i,j \leq n; 1 \leq \ell \leq s]$$
be the commutative algebra of polynomial functions, with coefficients
in $k$, on $\M$.

(iii) For $1 \leq \ell \leq s$, let $\x_\ell$ denote the generic matrix
$$\big(x_{ij}(\ell)\big) \in M_n(B).$$

(iv) Let $K$ be a commutative ring, and let $K\{X_1,\ldots,X_m\}$ be the free
associative $K$-algebra in the noncommuting variables $X_1,\ldots,X_m$. Given
a $K$-module $M$, we will also regard $K\{X_1,\ldots,X_m\}$ as an algebra of
noncommutative polynomial functions from $M^m$ to $M$.

(v) Choose $f_1,\ldots,f_t \in k\{ X_1, \ldots , X_s \}$, and set
$$R = k\{ X_1, \ldots , X_s \}/ \langle f_1, \ldots, f_t \rangle .$$
We will let $X_1,\ldots,X_s$ also denote their images in $R$. 

Our aim now is to algorithmically determine whether irreducible
representations $R \rightarrow M_n(\k)$ exist, and to construct one if they
do.

\subhead 2.1 \endsubhead (i) Each representation $\rho \colon R
\rightarrow M_n(\k)$ is determined exactly by the point
$$(\rho(X_1),\ldots,\rho(X_s)) \in \M.$$
In particular, the set of representations $R \rightarrow M_n(\k)$ can
be identified with
$$\left. \left\{  (\Gamma_1,\ldots,\Gamma_s) \in \M \; \right|
f_1(\Gamma_1,\ldots,\Gamma_s) = \cdots =
f_t(\Gamma_1,\ldots,\Gamma_s) = 0  \right\},$$
which is equal to the closed subvariety $V(\Rel(B))$ of $\M$, where
$\Rel(B)$ is the ideal of $B$ generated by the entries of the matrices
$$f_1(\x_1,\ldots,\x_s),\ldots,f_t(\x_1,\ldots,\x_s)$$
in $M_n(B)$.

(ii) Let $P$ denote the set of $s$-tuples $(\Gamma_1, \ldots,
\Gamma_s) \in \M$ for which the $\k$-algebra generated by the
$\Gamma_1,\ldots,\Gamma_s$ is not equal to $M_n(\k)$. Since $P$ is
equal to the set of $s$-tuples of simultaneously
block-upper-triangularizable matrices (for non-$\nxn$ blocks), we see
that $P$ is a closed subvariety of $\M$.

(iii) Suppose that $P$ is defined by the equations $g_1 = \cdots = g_q
= 0$ in $B$. By definition, there exists an irreducible representation
$R \rightarrow M_n(\k)$ if and only if $V(\Rel(B)) \not\subseteq P$.
Therefore, there exists an irreducible representation $R \rightarrow
M_n(\k)$ if and only if at least one $g_i$ is not contained in the
radical of $\Rel(B)$. Consequently, the radical membership algorithm
can be used to determine whether or not there exists an irreducible
$n$-dimensional representation of $R$.

It remains, then, to specify a set of defining equations for $P$.

\subhead 2.2 \endsubhead Let $K$ be a field, and let $A$ be a
$K$-subalgebra of $M_n(K)$. Suppose further that $A$ is generated, as
a $K$-algebra, by the set $G$. Let $p = n^2$. It is easy to see that
$A$ is $K$-linearly spanned by
$$\left\{ a_1 \cdots a_i \; \mid \; a_1,\ldots,a_i \in G, \; 0 \leq i
< p \right\},$$
where the product corresponding to $i=0$ is the identity matrix. It
follows from \cite{\Pap} that the preceding conclusion remains true if
we instead use
$$p = n\sqrt{2n^2/(n-1)+ 1/4} + n/2 -2 .$$
(Moreover, by the Cayley-Hamilton Theorem, we can always replace
$a^n$, for $a \in A$, by a $K$-linear combination of $1, a, a^2,
\ldots , a^{n-1}$.)

\subhead 2.3 \endsubhead We now turn to polynomial identities. Our
brief treatment here is distilled from \cite{\McCRob, Chapter 13}
(cf\. \cite{\Fortwo, \Row}). Let $A$ be a $k$-algebra, and let $g \in
\Z \{Y_1,\ldots, Y_m \}$. 

(i) If $X$ is a subset of $A$ then the set $\{ g(a_1,\ldots,a_m) \mid
a_1,\ldots,a_m \in X \}$ will be designated $g(X)$. 

(ii) The {\sl $m$th standard identity\/} is
$$s_m \; = \; \sum _{\sigma \in S_m} (\sgn
\sigma)Y_{\sigma(1)}\cdots Y_{\sigma(m)} \; \in \;
\Z\{Y_1,\ldots,Y_m\} .$$
Observe that $s_m\colon A^m \rightarrow A$ is $Z(A)$-multilinear and
alternating, where $Z(A)$ denotes the center of $A$.

\subhead 2.4 \endsubhead (i) Let $K$ be a commutative ring. The
Amitsur-Levitzky Theorem (see, e.g., \cite{\McCRob, 13.3.3}) ensures
that
$$s_{2m'}(M_m(K)) = 0$$
for all $m' \geq m$. Moreover, $s_{2m'}(M_m(K)) \ne 0$ for all $m' <
m$ (e.g., \cite{\McCRob, 13.3.2}).

(ii) Let $K$ be a field, and suppose that $A$ is a proper
$K$-subalgebra of $M_n(K)$. Let $J$ denote the Jacobson radical of
$A$. The semisimple algebra $A/J$ will embed (as a non-unital subring)
into a direct sum of copies of $M_m(K)$, for some $m < n$. It
therefore follows from (i) that $s_{2(n-1)}(A/J) = 0$, and so
$s_{2(n-1)}(A) \subseteq J$.

\subhead 2.5 \endsubhead Let $K$ be a field, and let $A$ be a
$K$-subalgebra of $M_n(K)$. Let $J$ denote the Jacobson radical of
$A$, and set
$$L = A \cdot s_{2(n-1)}(A), $$
a left ideal of $A$.

(i) Suppose that $A$ is a proper subalgebra of $M_n(K)$. Then $L$ is a
left ideal contained within $J$, by (2.4), and hence $L$ is a
nilpotent left ideal of $A$. In particular, every matrix in $L$ has
trace zero.

(ii) If $A = M_n(K)$ then $L$ is a left ideal of $M_n(K)$, and at
least one matrix in $L$ has nonzero trace.

(iii) Let $T = \{b_1,\ldots,b_N\}$ be a $K$-linear spanning set for
$A$. Set
$$V = \left\{ b_{m_0} \cdot s_{2(n-1)}(b_{m_1}, \ldots, b_{m_{2(n-1)}}) \;
\mid \; 1 \leq m_0 \leq N, \quad 1 \leq m_1 < \cdots < m_{2(n-1)} \leq N
\right\},$$
and note that $V$ is a linear generating set for $L$. (Recall that
$s_{2(n-1)}$ is multilinear and alternating.)

(iv) We conclude that $A$ is a proper subalgebra of $M_n(K)$ if and
only if $\{\trace(v) \mid v \in V\} = \{0\} $.

(v) Suppose that $A$ is generated, as a $K$-algebra, by $G \subseteq
M_n(K)$. Choosing $p$ as in (2.2), we may take
$$T = \left\{ a_1 \cdots a_i \; \mid \; a_1,\ldots,a_i \in G, \; 0
\leq i < p \right\}.$$

\subhead 2.6 \endsubhead We now prove that the proposed algorithm
satisfies the claims made in (1.1). A summary of the procedure, and
comments on it, will be presented in the next section.

(i) Let $p$ be as in (2.2), and set
$$S = \left\{\x_{\ell_1}\cdots \x_{\ell_i} \; \mid \; 0 \leq i < p  \right\}
\subseteq M_n(B).$$

(ii) Write $S = \{M_1,\ldots,M_N\}$, and set $U =$
$$\left\{ M_{m_0} \cdot s_{2(n-1)}(M_{m_1}, \ldots, M_{m_{2(n-1)}})
\; \mid \; 1 \leq m_0 \leq N, \quad 1 \leq m_1 < \cdots <
m_{2(n-1)} \leq N \right\}.$$

(iii) Recall $P \subseteq \M$ from (2.1ii). It follows from (2.5) that 
$$\{ \trace(u) = 0 \; | \; u \in U  \}$$
is a set of defining equations, in $B$, for $P$. 

(iv) Following (2.1), there exists an irreducible representation $R
\rightarrow M_n(\k)$ if and only $\trace(U) = \{ \trace(u) \; | \;
u \in U \}$ is not contained in the radical of $\Rel(B)$, and we
may therefore use the radical membership test to determine whether or
not $R$ has an irreducible $n$-dimensional representation.

(v) Suppose that $y \in \trace(U) \subseteq B$ is not contained in the
radical of $\Rel(B)$. Further suppose that $k[x]$ is equipped with a
factoring algorithm. Elimination methods can now be applied to find a
homomorphism $\varphi\colon B \rightarrow \k$ such that $y \not\in
\ker\varphi$ and such that $\Rel(B) \subseteq \ker\varphi$.  The
assignment
$$X_\ell \longmapsto \bigg(\varphi \big(x_{ij}(\ell)\big)\bigg) \in
M_n(\k),$$
for $1 \leq \ell \leq s$, will then produce an irreducible
$n$-dimensional representation of $R$.

(vi) Other sets of polynomials can be used to define $P$. For example, we can
rewrite the matrices in $S$ as $n^2{\times}1$ column matrices, and then
concatenate all possible combinations of $n^2$ of them, to form
$n^2{\times}n^2$-matrices over $B$. Letting $D$ denote the set of determinants
of these matrices, we see that $P = V(D)$. (My thanks to Zinovy Reichstein for
this observation.) The variety $P$ can also be defined using the well-known
central polynomials described, for example, in \cite{\Forone, \Raz}.

(vii) Suppose that $s = n = 2$. Note that $\Gamma_1,\Gamma_2 \in
M_2(\k)$ generate $M_2(\k)$ as a $\k$-algebra if and only if
$\Gamma_1$ and $\Gamma_2$ are not simultaneously upper
triangularizable. By considering the possible Jordan canonical forms
of $\Gamma_1$, it is not hard to verify that $\Gamma_1$ and $\Gamma_2$
generate $M_2(\k)$ if and only if $\det(\Gamma_1\Gamma_2 -
\Gamma_2\Gamma_1) \ne 0$. Therefore, in this case, $R$ has an
irreducible $2$-dimensional representation if and only if
$\det(\x_1\x_2 - \x_2\x_1)$ is not contained in the radical of
$\Rel(B)$. The reader is referred to \cite{\BouBou, \Fri} for a
complete discussion of similarity classes of $2{\times}2$ matrices.

\head 3. The Procedure \endhead 

\subhead 3.1 \endsubhead We now outline a procedure, based on the
preceding section, that satisfies (1.1). A proof that the process
works follows from (2.6).

\subsubhead 1. Input \endsubsubhead (i) $n$ is a positive integer.

(ii) $k$ is a computable field, and $\k$ is the algebraic closure of $k$.

(iii) $R = k\{X_1,\ldots,X_s\}/\langle f_1,\ldots,f_t \rangle$.

\subsubhead 2. Notation \endsubsubhead (i) $B$ is the polynomial ring
in commuting variables $x_{ij}(\ell)$, for $1 \leq i,j \leq n$ and $1
\leq \ell \leq s$.

(ii) $M_n(B)$ is the $k$-algebra of $\nxn$ matrices over $B$, and
$\x_\ell$ denotes the $\ell$th generic matrix $(x_{ij}(\ell)) \in M_n(B)$.

(iii) $\Rel(B)$ denotes the ideal of $B$ generated by the entries of
$f_1(\x_1, \ldots, \x_s)$, $\ldots$, $f_t(\x_1,\ldots,x_s)$.

(iv) $s_{2(n-1)} = \sum _{\sigma \in S_{2(n-1)}} (\sgn
\sigma)Y_{\sigma(1)}\cdots Y_{\sigma(2(n-1))} \in
\Z\{Y_1,\ldots,Y_{2(n-1)}\}$, the $2(n-1)$th standard polynomial.

\subsubhead 3. Decision \endsubsubhead (i) For $p = 4$ when $n=2$, and
for (e.g.) $p = n\sqrt{2n^2/(n-1)+ 1/4} + n/2 -2$ otherwise (see
(2.2)), set
$$S = \left\{ \x_{\ell_1}\cdots \cdots \x_{\ell_m} \; \mid \; m < p
\right\}.$$
(By the Cayley-Hamilton Theorem, we may -- for example -- exclude from
the preceding set those terms containing $\x_\ell^n$, for $1 \leq \ell
\leq s$.) Choose an ordering for $S$, say $S = \{ M_1,\ldots, M_N \}$.

(ii) Set $U=$
$$\left\{ M_{m_0} \cdot s_{2(n-1)}(M_{m_1}, \ldots, M_{m_{2(n-1)}})
\; \mid \; 1 \leq {m_0} \leq N, \quad 1 \leq m_1 < \cdots < m_{2(n-1)}
\leq N \right\}.$$
(Recall that $s_{2(n-1)}$ is alternating.)

(iii) Applying the radical membership algorithm, determine whether any
elements in $\trace (U)$ are contained in the radical of
$\Rel(B)$. (Not every element of $U$ needs its trace evaluated, since
$\trace(YZ) = \trace(ZY)$). Also, for $y \in \trace(U)$, it may be
easier to test whether the image of $y$ in $B/\Rel(B)$ is contained in
the nilradical of $B/\Rel(B)$; working modulo $\Rel(B)$, the generic
matrix arithmetic can often be significantly simplified.) If every
element in $\trace(U)$ is contained in the radical of $\Rel(B)$ then
there exist no irreducible representations $R \rightarrow M_n(\k)$;
see (2.6iv). If at least one element in $\trace(U)$ is not contained
in the radical of $\Rel(B)$, then there exist irreducible
representations $R \rightarrow M_n(\k)$, and we may proceed to step 4.

\subsubhead 4. Construction \endsubsubhead If $k[x]$ is equipped with
a factoring algorithm, and if $y \in \trace(U)$ is not contained in
the radical of $\Rel(B)$:

(i) Apply elimination methods to solve the $sn^2 + 1$ commutative polynomial
equations, in $B[z]$, obtained by setting $yz-1$ and the entries of
$f_1(\x_1,\ldots,\x_s),\ldots,f_t(\x_1,\ldots,\x_s)$ equal to zero. In this
solution, say, $x_{ij}(\ell) = \lambda_{ij}(\ell) \in \k$, for $1 \leq i,j
\leq n$ and $1 \leq \ell \leq s$.

(ii) The representation
$$\matrix R & \longrightarrow & M_n(\k) \\ \\ X_\ell & \longmapsto &
(\lambda_{ij}(\ell)) \endmatrix $$
is irreducible, by (2.6v).

\subhead 3.2 Remarks \endsubhead (i) It is sensible, in practice, to first
look for irreducible representations $\rho\colon R \rightarrow M_n( k)$ under
simplifying assumptions. For example, one can initially suppose that one (or
more) of the $\rho(X_\ell)$ are diagonal, or that a subset of the images of
the $\rho(X_\ell)$ are triangular; see example (4.2). (Of course, for
any commuting subset of the generators $X_\ell \in R$, there is no loss of
generality in assuming that the images are all upper triangular.)

(ii) Roughly speaking, the cost of employing this procedure depends on
the degrees of the polynomials involved in applications of the radical
membership algorithm. Note, in general, that the polynomials in
$\trace(U)$ may have degree $p^{2n-1}$. Another consideration will be
the number of polynomials in $\trace(U)$ to which the radical
membership algorithm, modulo $\Rel(B)$, is actually applied. This
quantity appears difficult -- in general -- to precisely estimate and
can vary greatly for different choices of $f_1,\ldots,f_t$; see
example (4.3). Observe, if the number of elements of $S$ used in step
3 is equal to $q$, that the number of terms $M_{m_0} \cdot
s_{2(n-1)}(M_{m_1}, \ldots, M_{m_{2(n-1)}})$ is $q{\binom q
{2(n-1)}}$.

(iii) Recalling (2.6vi), one can use $D$ instead of $\trace(U)$ in
steps 3 and 4. In general, the polynomials in $D$ can have degree
$p^{n^2}$, and if $q$ is the number of elements from $S$ used in this
approach then there will be ${\binom q {n^2}}$ polynomials to which
the radical membership algorithm must be applied. 

(iv) By (2.6vii), when $s = n = 2$, we can replace $\trace(U)$ with
the single polynomial $\det(\x_1\x_2 - \x_2\x_1)$.

\head 4. Examples \endhead

Retain the notation of the previous sections. 

\subhead 4.1 \endsubhead We begin with $2$-dimensional
representations.         

(i)  Set 
$$\x_1 = \x = \bmatrix x_{11} & x_{12} \\ x_{21} & x_{22} \endbmatrix,
\quad \x_2 = \y =\bmatrix y_{11} & y_{12} \\ y_{21} & y_{22}
\endbmatrix, \quad \x_3 = \z =\bmatrix z_{11} & z_{12} \\ z_{21} & z_{22}
\endbmatrix. $$
and 
$$B = \Q [ x_{11},x_{12},x_{21},x_{22}, y_{11},y_{12},y_{21},y_{22},
z_{11},z_{12},z_{21},z_{22}].$$

(ii) The value of $p$, as defined in (3.1.3i), is $4$, and
$$S = \left\{  abc \; | \; a,b,c \in \{ 1, \x, \y, \z\} \right\}.$$

(iii) Here, $s_{2(n-1)} = s_2(a,b) = ab-ba$ is the commutator. As in
(3.1.3ii), order $S$, and let
$$U = \{ a(a'a''-a''a') \; \mid \; a,a',a'' \in S, \quad a' < a'' \}
.$$

(iv) Now set
$$R = \Q\{X,Y,Z\}/\langle (XY - YX)^2, (XZ - ZX)^2, (YZ - ZY)^2
\rangle .$$
Using Macaulay2, we found that
$$T=\trace(\x(\y\z-\z\y))$$
is not contained in the radical of $\Rel(B)$; see
(3.1.3iii). Therefore, $R$ has an irreducible $2$-dimensional
representation. For example,
$$X \mapsto \bmatrix 2 & \hfill 0 \\ 0 & -2 \endbmatrix, \quad Y \mapsto
\bmatrix 2 & \hfill 2 \\ 0 & -2 \endbmatrix, \quad Z \mapsto \bmatrix 1 & 0 \\
2 & 2 \endbmatrix,$$
defines an irreducible $2$-dimensional representation in which $T \mapsto
16$. By (3.2iv
), the subalgebras of $R$ generated by any two of the
generators $X, Y, Z$ have no irreducible $2$-dimensional representations.

\subhead 4.2 \endsubhead Continue to let $R$ be as in (4.1); we now
consider the case when $n = 3$. Let $\x$, $\y$, and $\z$ denote the
$3{\times}3$ generic matrices respectively corresponding to $X$, $Y$,
and $Z$. 

To make the calculations more manageable, one can first check to see
if $R$ has an irreducible $3$-dimensional represention in which $X$ is
diagonal, $Y$ is upper triangular, and $Z$ is lower triangular. With
this simplification, using Macaulay2, we found that
$$T = \trace(\x\cdot s_4 (\y, \z, \x\y, \x\z ))$$
is not contained in the radical of $\Rel(B)$, and so $R$ must have a
$3$-dimensional irreducible representation. For instance,
$$X \mapsto \bmatrix 2 & \hfill 0 & 0 \\ 0 & -2 & 0 \\ 0 & \hfill 0 & 2
\endbmatrix, \quad Y \mapsto \bmatrix 2 & -1 & 2 \\ 0 & -2 & 8 \\ 0 &
\hfill 0 & 2 \endbmatrix, \quad Z \mapsto \bmatrix 2 & 0 & \hfill 0 \\
4 & 2 & \hfill 0 \\ 2 & 2 & -1 \endbmatrix$$
produces a $3$-dimensional irreducible representation in which $T
\mapsto 8192$.

\subhead 4.3 \endsubhead Set $n =3$ and $R =$
$$\Q\{a,b,X,Y\} \left/ \left\langle \matrix X^2 - a , Y^2 - b
, \\ \text{$uv - vu$ for $u \in \{ a,b\}$ and $v \in
\{a,b,X,Y$\}} \endmatrix \right\rangle\right. .$$
Let $\x$ and $\y$ be the $3{\times}3$ generic matrices corresponding,
respectively, to $X$ and $Y$.

Following (3.1.3i), we can take $8 < p < 9$. In view of the defining
relations for $R$, we may now set $S = \{ M_1, \ldots , M_{17} \} =$
$$\left\{ \matrix 1, \x, \y, \x\y, \y\x, \x\y\x, \y\x\y, \x\y\x\y,
\y\x\y\x, \x\y\x\y\x, \y\x\y\x\y, \x\y\x\y\x\y, \\ \y\x\y\x\y\x,
\x\y\x\y\x\y\x, \y\x\y\x\y\x\y, \y\x\y\x\y\x\y\x, \x\y\x\y\x\y\x\y
\endmatrix \right\},$$
as in (3.1.3i). Following (3.1.3ii),
$$U = \{ M_{m_0} \cdot s_4(M_{m_1},M_{m_2},M_{m_3},M_{m_4}) \; \mid \;
1 \leq m_0 \leq 17, \; 1 \leq m_1 < m_2 < m_3 < m_4 \leq 17 \}.$$
Using Macaulay2, we checked directly that every member of $\trace(U)$
is contained in the radical of $\Rel(B)$. Therefore, by (3.1.3iii),
there exist no $3$-dimensional irreducible representations of $R$.

\enddocument